\documentclass[12pt]{article}

\usepackage{latexsym}
\newcommand{\Z}{\ensuremath{\mathsf{Z}}}
\newcommand{\R}{\ensuremath{\mathsf{R}}}
\newcommand{\C}{\ensuremath{\mathsf{C}}}
\newcommand{\N}{\ensuremath{\mathsf{N}}}
\newtheorem{thm}{Theorem}[section]

\begin{document}
\begin{center}
\Large{Strong perforation in infinitely generated $\mathrm{K}_0$-groups of
  simple $C^*$-algebras}

\vspace{5mm}

\large{Andrew S. Toms}
\normalsize

\end{center}
\begin{abstract}

Let $(G,G^+)$ be an ordered abelian group.  We say that $G$ has strong
perforation if there exists $x \in G$, $ x \notin G^+$, such that $nx
\in G^+$, $nx \neq 0$ for some natural number $n$.  Otherwise, the
group is said to be weakly unperforated.  Examples of simple
$C^*$-algebras whose ordered $\mathrm{K}_0$-groups have this property and for which the
entire order structure on $\mathrm{K}_0$ is known have, until now, been restricted
to the case where $\mathrm{K}_0$ is group isomorphic to the integers.
We construct simple, separable, unital $C^*$-algebras with
strongly perforated $\mathrm{K}_0$-groups group isomorphic to an arbitrary
infinitely generated subgroup of the rationals, and determine the
order structure on $\mathrm{K}_0$ in each case.

\end{abstract}

\section{Introduction}

Elliott's classification of $\mathrm{AF}$ $C^*$-algebras via the 
$\mathrm{K}_0$-group ([2]) began a widespread effort to classify
nuclear $C^*$-algebras.  The $\mathrm{K}_0$-group, which is
an ordered group for stably finite $C^*$-algebras ([1]), has figured
prominently in almost all work on this problem.  (For an overview of 
the classification problem for nuclear $C^*$-algebras,
see [3].)  So far, every result on the classification of $C^*$-algebras has
required the assumption that the ordered $\mathrm{K}_0$-group be 
weakly unperforated whenever it is not zero.  This assumption was shown to be non-trivial
by Villadsen ([8]);  the ordered abelian group $\mathsf{Z}_n := 
( \mathsf{Z}, \{0,n,n+1,\ldots\} )$ may arise as a saturated
sub-ordered group of the $\mathrm{K}_0$-group of a simple nuclear
$C^*$-algebra.  In [4], Elliott and Villadsen refined the results of
[8] to obtain, for each natural number $n$, a simple nuclear $C^*$-algebra
$A_n$ whose ordered $\mathrm{K}_0$-group is order isomorphic to $\mathsf{Z}_n$.
This result was further generalised by the author in [7], where it was 
shown that a certain class of order structures on the integers (which might
possibly comprise all such order structures giving a simple ordered group)
could arise as the ordered $\mathrm{K}_0$-group of a simple nuclear
$C^*$-algebra.

The classification of a category by an invariant is not
complete until one knows the range of the invariant, and any 
classification of simple nuclear stably finite $C^*$-algebras 
will necessarily capture the ordered $\mathrm{K}_0$-group.  Thus,
the range of the $\mathrm{K}_0$ functor bears investigation. 
This range is known when $\mathrm{K}_0$
is a weakly unperforated ordered group, whence our interest
in instances of the ordered 
$\mathrm{K}_0$-group which exhibit strong perforation.

\emph{Acknowledgements}:  The author would like to acknowledge the 
support of NSERC (PGS A, PGS B, Postdoctoral Fellowship) and the
Israel Halperin Graduate Award (University of Toronto).

\section{Essential Results}

In this section we review results from [4] that will be used in the sequel.

Let $C$, $D$ be $C^{*}$-algebras, and let $\phi_{0}$, $\phi_{1}$ be
$*$-homomorphisms from $C$ to $D$.  The generalised mapping torus of
$C$ and $D$ with respect to $\phi_{0}$ and $\phi_{1}$ is 
\begin{displaymath}
A:=\{(c,d)|d \in C([0,1];D), \ c \in C, \ d(0)=\phi_{0}(c), \ d(1)=\phi_{1}(c)\}
\end{displaymath}
We will write $A(C,D,\phi_{0}, \phi_{1})$ for $A$ when clarity demands
it.  We now list without proof some theorems, specialised to our needs, 
which will be used in the sequel.
 
\begin{thm}[{Elliott and Villadsen ([4]), Sec. 2, Thm. 2}]
The index map $b_{*} : \mathrm{K}_{*}C \rightarrow 
\mathrm{K}_{1-*} \mathrm{S}D = \mathrm{K}_{*}D$ in the six term 
periodic sequence for the extension
\begin{displaymath}
0 \rightarrow \mathrm{S}D \rightarrow A \rightarrow C \rightarrow 0
\end{displaymath}
is the difference
\begin{displaymath}
\mathrm{K}_{*} \phi_{1} - \mathrm{K}_{*} \phi_{0} : \mathrm{K}_{*} C \rightarrow \mathrm{K}_{*} D.
\end{displaymath}
Thus, the six-term exact sequence may be written as the short exact sequence
\begin{displaymath}
0 \rightarrow \mathrm{Coker} b_{1-*} \rightarrow
\mathrm{K}_* A \rightarrow \mathrm{Ker} b_{*} \rightarrow 0.
\end{displaymath}
In particular, if $b_{1-*}$ is surjective, then $\mathrm{K}_{*} A$ 
is isomorphic to its image, $\mathrm{Ker} b_{*}$, in $\mathrm{K}_{*} C$.
 
Suppose that cancellation holds for each pair of 
projections in $D \otimes \mathcal{K}$ obtained as the images under the maps $\phi_{0}$ and
$\phi_{1}$ of a single projection in $C \otimes \mathcal{K}$.
Then, if $b_{1}$ is 
surjective, 
\begin{displaymath}
{(    \mathrm{K}_{0} A)}^{+} \cong {(    \mathrm{K}_{0} C)}^{+} \cap     \mathrm{K}_0(e_{\infty})(\mathrm{K}_{0} A),
\end{displaymath}
where $e_{\infty}$ denotes the evaluation of $A$ at the fibre at infinity.
\end{thm}

\begin{thm}[{Elliott and Villadsen ([4]), Sec. 3, Thm. 3}]
Let $A_{1}$ and $A_{2}$ be building block algebras as described above,

\begin{displaymath} A_{i}=A(C,D,\phi_{0}^{i},\phi_{1}^{i}),\ \ i=1,2.
\end{displaymath}

Let there be given three maps between the fibres,
\begin{displaymath}
\begin{array}{rrr}
\gamma : & C_{1} \rightarrow C_{2}, &  \\
\delta, \delta ' : & D_{1} \rightarrow D_{2}, & 
\end{array}
\end{displaymath}
such that $\delta$ and $\delta '$ have mutually orthogonal images, and
\begin{displaymath} \delta \phi_{0}^{1} + \delta ' \phi_{1}^{1} = \phi_{0}^{2} \gamma,
\end{displaymath}
\begin{displaymath} \delta \phi_{1}^{1} + \delta ' \phi_{0}^{1} = \phi_{1}^{2} \gamma.
\end{displaymath}

Then there exists a unique map 
\begin{displaymath} \theta : A_{1} \rightarrow A_{2},
\end{displaymath}
respecting the canonical ideals, giving rise to the map $\gamma : C_{1} \rightarrow
C_{2}$ between the quotients (or fibres at infinity), and such that for any $0 < s< 1$,
if $e_{s}$ denotes evaluation at $s$,
\begin{displaymath} e_{s} \theta = \delta e_{s} + \delta ' e_{1-s}.
\end{displaymath}

\end{thm}

Let $A_{1}$ and $A_{2}$ be building block algebras as in Theorem 2.1
with $\theta : A_{1} \rightarrow A_{2}$ as in Theorem 2.2.
Let there be given a map $\beta : D_{1} \rightarrow C_{2}$ such that the composed map
$\beta \phi_{1}^{1}$ is a direct summand of the map $\gamma$, and such that the 
composed maps $\phi_{0}^{2} \beta$ and $\phi_{1}^{2} \beta$ are direct summands of
the maps $\delta '$ and $\delta$, respectively.  Suppose that the decomposition of
$\gamma$ as the orthogonal sum of $\beta \phi_{1}^{1}$ and another map is such that
the image of the second map is orthogonal to the image of $\beta$.  (Note that this
requirement is automatically satisfied if $C_{1}$, $D_{1}$, and the map $\beta \phi_{1}^{1}$
are unital.)

Let
\begin{displaymath} A_{1} \stackrel{\theta_{1}}{\rightarrow} A_{2} \stackrel{\theta_{2}}{\rightarrow} \cdots
\end{displaymath}
be a sequence of separable building block $C^{*}$-algebras,
\begin{displaymath} A_{i} = A(C_{i},D_{i},\phi_{0}^{i}, \phi_{1}^{i}), \ i=1,2, \ldots
\end{displaymath} 
with each map $\theta_{i} : A_{i} \rightarrow A_{i+1}$ obtained by the
construction of Theorem 2.2.  For each $i=1,2, \ldots$ let $\beta_{i} 
: D_{i} \rightarrow C_{i+1}$ be a map verifying the hypotheses of the
preceding paragraph.

Suppose that for every $i=1,2, \ldots$, the intersection of the kernels of the boundary maps $\phi_{0}^{i}$ and 
$\phi_{1}^{i}$ from $C_{i}$ to $D_{i}$ is zero.

Suppose that, for each $i$, the image of each of $\phi_{0}^{i+1}$ and $\phi_{1}^{i+1}$ generates $D_{i+1}$ as a
closed two-sided ideal, and that this is in fact true for the restriction of $\phi_{0}^{i+1}$ and 
$\phi_{1}^{i+1}$ to the smallest direct summand of $C_{i+1}$ containing the image of $\beta_{i}$.  Suppose
that the closed two-sided ideal of $C_{i+1}$ generated by the image of $\beta_{i}$ is a direct summand. 

Suppose that, for each $i$, the maps $\delta_{i}^{'} - \phi_{0}^{i} \beta_{i}$ and 
$\delta_{i} - \phi_{1}^{i} \beta_{i}$ from $D_{i}$ to $D_{i+1}$ are injective.  

Suppose that, for each $i$, the map $\gamma_{i} - \beta_{i} \phi_{1}^{i}$ takes each non-zero direct summand
of $C_{i}$ into a subalgebra of $C_{i+1}$ not contained in any proper closed two-sided ideal. 

Suppose that, for each $i$, the map $\beta_{i} : D_{i} \rightarrow C_{i+1}$ can be deformed---inside the 
hereditary sub-$C^{*}$-algebra generated by its image---to a map $\alpha_{i} : D_{i} \rightarrow C_{i+1}$
with the following property:  There is a direct summand of $\alpha_{i}$, say ${\bar{\alpha}}_{i}$, such
that ${\bar{\alpha}}_{i}$ is non-zero on an arbitrary given element $x_{i}$ of $D_{i}$, and has image a 
simple sub-$C^{*}$-algebra of $C_{i+1}$, the closed two-sided ideal generated by which contains the image 
of $\beta_{i}$.

\begin{thm}[{Elliott and Villadsen ([4]), Sec. 5, Thm 5}]
If the hypotheses above are satisfied, there is a map $\theta_{i}^{'}$ homotopic inside $A_i$ to $\theta_{i}$
for each $i$ such that the inductive limit of the sequence
\begin{displaymath} A_{1} \stackrel{\theta_{1}^{'}}{\rightarrow} A_{2} \stackrel{\theta_{2}^{'}}{\rightarrow}
\cdots
\end{displaymath}
is simple.
\end{thm}

\section{Infinitely Generated Subgroups of the Rational Numbers}

 A generalised integer is a symbol 
$\mathbf{n} = a_{1}^{n_{1}} a_{2}^{n_{2}} a_{3}^{n_{3}} \ldots$, where the 
$a_{i}$'s are pairwise distinct prime numbers and each $n_{i}$ is
either a non-negative integer or $\infty$.  
The subgroup $G_{\mathbf{n}}$ of the rational numbers associated to the generalised
integer $\mathbf{n}$ is the group of all rationals whose denominators 
(when in lowest terms) are products of powers of the $a_{i}$'s not 
exceeding $a^{n_{i}}$.  If $n_{i} = \infty$, then an arbitrarily large power of $a_{i}$ 
may appear in the denominator.  

\begin{thm} For each pair $(\mathbf{n},k)$ consisting of a 
generalised integer $\mathbf{n}$ and a positive rational 
$k < 1$, there exists a simple, separable, unital, nuclear 
$C^{*}$-algebra $A_{(\mathbf{n},k)}$ such that
\begin{displaymath}
(    \mathrm{K}_{0}(A_{(\mathbf{n},k)}),
{\mathrm{K}_{0}(A_{(\mathbf{n},k)})}^{+}, [1_{A_{(\mathbf{n},k)}}] 
) = ( G_{\mathbf{n}}, G_{\mathbf{n}} \cap (k,\infty), 1 ).
\end{displaymath}
\end{thm}

\noindent
\textbf{Proof:}
Given a 2-tuple $(\mathbf{n},k)$ we will construct a sequence
\begin{displaymath}
A_{1} \stackrel{\theta_{1}}{\rightarrow} A_{2} \stackrel{\theta_{2}}{\rightarrow} \cdots
\end{displaymath}
where $A_{j} = A(C_{j},D_{j},\phi_{0}^{j},\phi_{1}^{j})$, and the
$\theta_{j}$ constructed as in Theorem 2.2 from maps 
\begin{displaymath}
\gamma_{j} : C_{j} \rightarrow C_{j+1}, \  \ \delta_{j}, \delta_{j}^{'} : D_{j} \rightarrow D_{j+1}.
\end{displaymath}
In order to obtain a simple inductive limit, we will require a map
\begin{displaymath}
\beta_{j} : D_{j} \rightarrow C_{j+1}
\end{displaymath}
having the properties listed in Section 2.  

For each $j$ let
\begin{displaymath}
C_{j} = p_{j} (\mathrm{C}(X_{j}) \otimes \mathcal{K}) p_{j}
\end{displaymath}
where $p_{j}$ is a projection in $\mathrm{C}(X_{j}) \otimes
\mathcal{K}$ and $\mathcal{K}$ denotes the compact operators.  Express $k$ in lowest terms, 
say $\frac{a}{b}$, and set $X_{1} = {\mathrm{S}^{2}}^{\times (a+1)}$. 
Let $X_{j+1} = {X_{j}}^{\times n_{j}}$, where $n_{j}$ is a natural number 
to be specified.  

Let $D_{j} = C_{j} \otimes \mathrm{M}_{\mathrm{dim}(p_{j})k_{j}}$, where $k_{j}$ is a natural number to be 
specified.    
Let $\mu_{j}$ and $\nu_{j}$ be maps
from $C_{j}$ to $C_{j} \otimes \mathrm{M}_{\mathrm{dim}(p_{j})}$ given
by 
\begin{displaymath}
\mu_{j}(a) = p_{j} \otimes a(x_{j}) \cdot 1_{\mathrm{dim}(p_{j})}
\end{displaymath}
(where $x_{j}$ is a point to be specified in $X_{j}$ and $1_{\mathrm{dim}(p_{j})}$ is the 
unit of $\mathrm{M}_{\mathrm{dim}(p_{j})}$) and
\begin{displaymath}
\nu_{j}(a) = a \otimes 1_{\mathrm{dim}(p_{j})}.
\end{displaymath}
For $t \in \{0,1\}$, let $\phi_{j}^{t}:C_j \rightarrow D_j$ be the 
direct sum of $l_{j}^{t}$ and $k_{j} - l_{j}^{t}$ 
copies of $\mu_{j}$ and $\nu_{j}$, respectively,
where the $l_{j}^{t}$ are non-negative integers such that $l_{j}^{0} \neq l_{j}^{1}$ for all $j \geq 1$.

Note that both $C_{j}$ and $D_{j}$ are unital, as are the maps $\phi_{j}^{t}$.  
The $\phi_{j}^{t}$ are also injective and as such satisfy the hypotheses 
of Section 2 concerning them alone.

By Theorem 2.1, for each $e \in     \mathrm{K}_{0}(C_{j})$,
\begin{displaymath}
\begin{array}{rl}
b_{0}(e) & = (l_{j}^{1} - l_{j}^{0})(    \mathrm{K}_{0}(\mu_{j}) -     \mathrm{K}_{0}(\nu_{j}))(e) \\
& = (l_{j}^{1} - l_{j}^{0})(\mathrm{dim}(p_{j}) \cdot     \mathrm{K}_{0}(p_{j}) - \mathrm{dim}(p_{j}) \cdot e).  
\end{array}
\end{displaymath}
Since $l_{j}^{1} - l_{j}^{0}$ is non-zero for every $j$ and $    \mathrm{K}_{0}(X_{j})$ is torsion free, 
$b_{0}(e) = 0$ implies that $e$ belongs to the maximal free cyclic subgroup of $    \mathrm{K}_{0}(C_{j})$ containing $    \mathrm{K}_{0}(p_{j})$.
As $\mathrm{K}_{1}(C_{j}) = 0$, $b_{1}$ is surjective.  $\mathrm{K}_{0}(A_{j})$
is thus group isomorphic (by Theorem 2.1) 
to its image, in $\mathrm{K}_{0}(C_{j})$ --- which is isomorphic as a group to \Z.

In order for $    \mathrm{K}_{0}(A_{j})$ to be isomorphic as an ordered group to its image 
in $    \mathrm{K}_{0}(C_{j})$, with the relative order, it is sufficient 
(by Theorem 2.1) that for any projection $q$ in $C_{j} \otimes \mathcal{K}$ such that the 
images of $q$ under $\phi_{j}^{0} \otimes \mathrm{id}$ and
$\phi_{j}^{1} \otimes \mathrm{id}$ have the same $    \mathrm{K}_{0}$ class, these 
images be in fact equivalent.  For any such $q$, the image of $    \mathrm{K}_{0}(q)$
under $b_{0} = \mathrm{K}_{0}(\phi_{j}^{1}) - \mathrm{K}_{0}(\phi_{j}^{0})$ is zero, so that 
$\mathrm{K}_{0}(q)$ belongs to $\mathrm{Ker} b_{0}$.  It will be clear
from the construction below that the dimension of both $\phi_{j}^{1}(q)$ and 
$\phi_{j}^{0}(q)$ is at least half the dimension of $X_{j}$.  Thus, by
Theorem 8.1.5 of [5], $\phi_{j}^{1}(q)$ and $\phi_{j}^{0}(q)$ are equivalent,
as they have the same $\mathrm{K}_{0}$ class.

Let us now specify the projection $p_{1}$.  Let $\xi$ be the Hopf line 
bundle over $\mathrm{S}^{2}$.  Set $g_{1} = [\xi^{\times a+1}] - [\theta_{a}]
\in \mathrm{K}^{0}(X_{1})$, where $[ \ \cdot \ ]$ denotes the stable 
isomorphism class of a vector bundle and $\theta_{l}$ denotes the 
trivial vector bundle of fiber dimension $l$.  By Theorem 8.1.5
of [5], we have that 
$(a+1) \cdot g_{1}$ and hence $b \cdot g_{1}$ are positive.  
Let $p_{1}$ be a projection in $\mathrm{C}(X_{1}) \otimes \mathcal{K}$ corresponding to the
$\mathrm{K}^{0}$ class 
$b \cdot g_{1}$.  By [8] we know that the ordered, saturated, free cyclic 
subgroup of $\mathrm{K}_{0}(C_{1})$ generated by $g_{1}$ is equal to
\begin{displaymath}
( \Z, \{0,a+1,a+2, \ldots \} ),
\end{displaymath}    
where the class of the unit is the integer $b \geq a+1$.

  Decompose
$b$ into powers of primes, $b = a_{i_{1}}^{m_{1}} a_{i_{2}}^{m_{2}} \ldots 
a_{i_{n}}^{m_{n}}$.  Set ${\mathbf{n}}^{'} = \frac{\mathbf{n}}{b}$, with the convention that 
$\infty - l = \infty$ for all natural numbers $l$.  Let $L_{j}$ be an 
enumeration of the primes appearing in ${\mathbf{n}}^{'}$ for $j \geq 2$,
$j \in \N$, and set $L_{1} = b$. 

We now define a family of continuous maps from $\mathrm{S}^{2}$ to $\mathrm{S}^{2}$, indexed 
by the integers, to be used in the construction of the maps
$\gamma_{j}$ from $C_{j}$ to $C_{j+1}$.  Consider $\mathrm{S}^{2}$ as being embedded 
in ${\R}^{3} = \C \times \R$ as the unit sphere with center the origin, with
the identification $(x,y,z) = ( x + yi, z)$.  
For each $\eta \in \N$, let $\omega_{\eta}^{'}: \C \times \R \longrightarrow \C \times \R$ 
be defined by $\omega_{\eta}^{'}(w,z) = (w^{\eta}/|w^{\eta -1}|,z)$ when $w \neq 0$ and otherwise
by $\omega_{\eta}^{'}(0,z) = (0,z)$. This defines a map from $\mathrm{S}^{2}$
to itself by restriction.  Let $\omega_{\eta}$ be the composition of $\omega_{\eta}^{'}$ with the antipodal map.
Note that $\omega_{\eta}^{'}$ is the suspension of the
$\eta^{\mathrm{th}}$ power map on $\mathrm{S}^{1}$, and thus
has the same degree, namely $- \eta$, as this map ([6]).  As 
the antipodal map has degree $-1$, the composed map
$\omega_{\eta}$ has degree $\eta$.  In the language of 
vector bundles, $\mathrm{K}^{0}(\omega_{\eta})([\xi]) = [\xi^{\otimes \eta}]$.

Define a map $\gamma_{j}^{'}$  
from $\mathrm{C}(X_{j})$ to $\mathrm{M}_{n_{j}} \otimes 
\mathrm{C}(X_{j+1}) = \mathrm{M}_{n_{j}}({\mathrm{C}(X_{j}}^{\otimes n_{j}})$ 
as follows:
\begin{displaymath}
\begin{array}{rl}
\gamma_{j}^{'}(f(x)) & = (f(\omega_{L_{j+1}}(x)) \otimes 1 \otimes \ldots \otimes 1)  \ \oplus \ 
                         (1 \otimes f(\omega_{L_{j+1}}(x)) \otimes \ldots \otimes 1) \oplus \ldots \\
                     & \ \ \ \ \ \ \ \ \  \ldots \oplus \ (1 \otimes 1 \otimes \ldots \otimes f(\omega_{L_{j+1}}(x))).
\end{array}
\end{displaymath}

Let
\begin{displaymath}
\beta_{j}^{'}  = 1 \cdot e_{x_{j}}
\end{displaymath}
be a map from $\mathrm{C}(X_{j})$ to $\mathrm{C}(X_{j+1})$, where $e_{x_{j}}$ denotes the evaluation 
of an element of $\mathrm{C}(X_{j})$ at a point $x_{j} \in X_{j}$ and $1$ is the unit of
$\mathrm{C}(X_{j+1})$.  Fix $x_1 \in \mathrm{S}^2$ and define
$x_{j+1}:=(\omega_{L_{j+1}}(x_j),\ldots,\omega_{L_{j+1}}(x_j)) \in
{X_j}^{\times n_j} = X_{j+1}$.

Let us define $\gamma_{j}: \mathrm{C}(X_{j}) \rightarrow
\mathrm{M}_{n_{j}}(\mathrm{C}(X_{j+1})) \otimes \mathrm{M}_{2}(\mathcal{K})$ 
inductively as the direct sum of two maps.  For the
first map, take the restriction to $C_{j} \subseteq \mathrm{C}(X_{j}) \otimes \mathcal{K}$ 
of the tensor product of $\gamma_{j}^{'}$ with the identity map from $\mathcal{K}$ to $\mathcal{K}$.  The second
map is obtained as follows:  compose the map $\phi_{j}^{1}$ with the direct sum 
of $q_{j}$ copies of the tensor product of $\beta_{j}^{'}$ with the identity map from
$\mathcal{K}$ to $\mathcal{K}$ (restricted to $D_{j} \subseteq \mathrm{C}(X_{j}) \otimes \mathcal{K}$), where $q_{j}$ 
is to be specified.  The induction consists of first considering the case $j=1$ (since
$p_{1}$ has already been chosen), then setting $p_{2} = \gamma_{j}(p_{1})$, so 
that $C_{2}$ is specified as the cut-down of $\mathrm{C}(X_{j}) \otimes \mathrm{M}_{2}(\mathcal{K})$, and continuing 
in this way.  

With $\beta_{j}$ taken to be the restriction to $D_{j} \subseteq \mathrm{C}(X_{j}) \otimes \mathcal{K}$ 
of $\beta_{j}^{'} \otimes \mathrm{id}$ we have, by construction, that $\beta_{j} 
\phi_{j}^{1}$ is a direct summand of $\gamma_{j}$ and, furthermore, the second 
direct summand and $\beta_{j}$ map into orthogonal blocks (and hence orthogonal 
subalgebras) as desired.  

We will now need to verify that $p_{j}$ has the following property:  the set of 
all rational multiples of $    \mathrm{K}_{0}(p_{j})$ in the ordered group $    \mathrm{K}_{0}(C_{j}) =
\mathrm{K}^{0}(X_{j})$ is isomorphic (as a sub ordered group) to 
\begin{displaymath}
( \Z, \{0, l_{j}+1, l_{j}+2,\ldots \} ),
\end{displaymath}
where
\begin{displaymath}
l_{j} := L_{j} l_{j-1}, \ l_{1} := a
\end{displaymath}
and the class of the unit (i.e., of $p_{j}$) is $\Pi_{k=1}^{j} L_{k}$. 

Our verification will proceed by induction.  The case $j=1$ has been established 
by construction.  Suppose that the assertion of the preceding paragraph holds for
all $p_{k}$, $k \leq j$.  Suppose further that the group of rational 
multiples of $\mathrm{K}_{0}(p_{k})$ (being isomorphic as a group to $\Z$) is generated by a
$\mathrm{K}_{0}$ class of the form $[\xi^{\times n}] - [\theta_{m}]$, where $m < n$ and  
(this is again true by construction for $k=1$).  We will show that $    \mathrm{K}_{0}(p_{j})$ 
has both the property of the preceding paragraph and the property just mentioned.

Let $g_{k} \in \mathrm{K}^{0}(X_{k})$ be the generator of the group of rational multiples 
of $p_{k}$.   
Note that, as is the case for all maps on $\mathrm{K}^{0}(\mathrm{S}^{2})$ 
induced by a continuous map from $\mathrm{S}^{2}$ to itself, $\mathrm{K}_{0}(\omega_{\eta})([\theta_{1}])
= [\theta_{1}]$.  Write $g_{k} = [\xi^{\times d_{k}}] - [\theta_{m_{k}}]$.  Then
\begin{displaymath}
    \mathrm{K}_{0}(\gamma_{j})(g_{j}) = [{(\xi^{\otimes L_{j+1}})}^{\times d_{j} n_{j}}] - [\theta_{m^{'}_{j+1}}]
\end{displaymath}
for some integers $d_{j} > 0$ and $m_{j+1}^{'}$.  We may assume that the multiplicity 
of the map $    \mathrm{K}_{0}(\gamma_{j})$  is divisible by $L_{j+1}$, as we have yet to specify $n_{j}$.  
We recall that for any integer $l$, the
$    \mathrm{K}_{0}$ class $[\xi^{\otimes l}]$ corresponds to the element $(1,l)$ in $\mathrm{K}^{0}(\mathrm{S}^{2}) 
= \langle [\theta_1] \rangle \oplus \langle e(\xi) \rangle$, which is also the difference of $\mathrm{K}_{0}$ classes
$l [\xi] - [\theta_{l - 1}]$.  Thus we have
\begin{displaymath}
    \mathrm{K}_{0}(\gamma_{j})(g_{j}) = L_{j+1} ([\xi^{\times (a+1) n_{1} n_{2} \ldots n_{j}}] - [\theta_{m_{j+1}}]).
\end{displaymath}
for some integer $m_{j+1}$.  Setting $g_{j} := [\xi^{\times (a + 1) n_{1} n_{2} 
\ldots n_{j}}] - [\theta_{m_{j+1}}]$, we have established that 
$\mathrm{K}_{0}(\gamma_{j})(g_{j}) = L_{j+1} g_{j+1}$ for all 
natural numbers $j$.

We now show that $n_{j}$ may be chosen so as to ensure that the maximal, 
free, cyclic subgroup of $    \mathrm{K}_{0} C_{j+1}$ generated by $g_{j+1}$ is indeed isomorphic 
as an ordered group to the integers with positive
cone $\{0, l_{j+1}+1, l_{j+1}+2, \ldots \}$.  That $\Pi_{k=1}^{j} L_{k}$ is the class 
of the unit follows directly from the fact that $L_{1} = b$ (the class of the unit in $    \mathrm{K}_{0} C_{1}$) and 
that $    \mathrm{K}_{0}(\gamma_{j})(g_{j}) = L_{j+1} g_{j+1}$.

As the Euler class of the Hopf line bundle on $\mathrm{S}^{2}$ is non-zero we have, by [8], that 
for $q, m, h \in \N$ such that $0 < h(q-m) < q$,
\begin{displaymath}
h ( [\xi^{\times q}] - [\theta_{m}] ) \notin {(\mathrm{K}^{0} {\mathrm{S}^{2}}^{\times q})}^{+}.
\end{displaymath}

To apply this we note that 
\begin{displaymath}
g_{j+1} = [\xi^{\times (a + 1) n_{1} n_{2} \ldots n_{j}}] - [\theta_{m_{j}}].
\end{displaymath}
With $q = (a+1) n_{1} n_{2} \ldots n_{j}$ and $m= m_{j}$ we wish to have
\begin{displaymath}
0 < l_{j}(q-m) < q
\end{displaymath}
as then $0 < h(q-m) < q$ for all $0 < h < l_{j}+1$.

Since
\begin{displaymath}
q-m = \mathrm{dim} g_{j+1} = \frac{n_{j} + k_{j} q_{j} \mathrm{dim} p_{j}}{L_{j+1}} \mathrm{dim} g_{j}
\end{displaymath}
we want
\begin{displaymath}
\mathrm{dim} g_{j+1} < \frac{(a+1) n_{1} n_{2} \ldots n_{j}}{l_{j+1}}.
\end{displaymath} 
Assume inductively that $n_{1}, n_{2}, \ldots, n_{j-1}$ have been chosen so that
\begin{displaymath}
\mathrm{dim} g_{j} < \frac{(a+1) n_{1} n_{2} \ldots n_{j-1}}{l_{j}}.
\end{displaymath}
Choose $n_{j}$ large enough so that
\begin{displaymath}
\frac{n_{j} + k_{j} q_{j} \mathrm{dim} p_{j}}{n_{j}} \mathrm{dim} g_{j} < \frac{(a+1) n_{1} n_{2} \ldots n_{j-1}}{l_{j}}.
\end{displaymath}
Then we have that
\begin{displaymath}
\frac{n_{j} + k_{j} q_{j} \mathrm{dim} p_{j}}{L_{j+1}} \mathrm{dim} g_{j} < \frac{(a+1) n_{1} n_{2} \ldots n_{j}}{L_{j+1} l_{j}}.
\end{displaymath}
Recalling that $l_{j+1} = L_{j+1} l_{j}$ we conclude that
\begin{displaymath}
\mathrm{dim} g_{j+1} = \frac{n_{j} + k_{j} q_{j} \mathrm{dim} p_{j}}{L_{j+1}} \mathrm{dim} g_{j} < 
\frac{(a+1) n_{1} n_{2} \ldots n_{j}}{l_{j+1}},
\end{displaymath}
as desired. 

Note that $\gamma_{j} - \beta_{j} \phi_{j}^{1}$ is non-zero and so, as required in the 
hypotheses of Theorem 2.4, takes $C_{j}$ into a subalgebra of 
$C_{j+1}$ not contained in any proper closed two-sided ideal.

It remains to construct maps $\delta_{j}$ and $\delta_{j}^{'}$ from $D_{j}$ 
to $D_{j+1}$ with orthogonal images such that
\begin{displaymath}
\delta_{j} \phi_{j}^{0} + \delta_{j}^{'} \phi_{j}^{1} = \phi_{j+1}^{0} \gamma_{j},
\end{displaymath}
\begin{displaymath}
\delta_{j} \phi_{j}^{1} + \delta_{j}^{'} \phi_{j}^{0} = \phi_{j+1}^{1} \gamma_{j},
\end{displaymath}
and $\phi_{j+1}^{0} \beta_{j}$ and $\phi_{j+1}^{1} \beta_{j}$ are direct summands of 
$\delta_{j}^{'}$ and $\delta_{j}$ respectively. To do this we shall have to modify 
$\phi_{j+1}^{0}$ and $\phi_{j+1}^{1}$ by inner automorphisms; this 
is permissible since it has no effect on $K$-groups.  The definition
of $\delta_{j}$ and $\delta_{j}^{'}$ along with the proof that they
satisfy the hypotheses of section 2 is taken from [4].

In order to carry out this step we define $x_{j+1} := \omega_{L_{j+1}}(x_{j})$, so that 
\begin{displaymath}
e_{x_{j+1}} \gamma_{j} = \mathrm{mult}(\gamma_{j}) e_{x_{j}},
\end{displaymath}
where $\mathrm{mult}(\gamma_{j})$ denotes the factor by which $\gamma_{j}$ multiplies 
dimension.  It follows that
\begin{displaymath}
\begin{array}{rl}
\mu_{j+1} \gamma_{j} & = p_{j+1} \otimes e_{x_{j+1}} \gamma_{j} \\
                     & = \gamma_{j}(p_{j}) \otimes \mathrm{mult}(\gamma_{j}) e_{x_{j}} \\
                     & = \mathrm{mult}(\gamma_{j}) \gamma_{j}(p_{j} \otimes e_{x_{j}}) \\
                     & = \mathrm{mult}(\gamma_{j}) \gamma_{j} \mu_{j},
\end{array}
\end{displaymath}
and
\begin{displaymath}
\begin{array}{rl}
\nu_{j+1} \gamma_{j} & = \gamma_{j} \otimes 1_{\mathrm{dim}( p_{j+1})} \\
                     & = \mathrm{mult}(\gamma_{j}) \gamma_{j} \otimes 1_{\mathrm{dim} (p_{j})} \\
                     & = \mathrm{mult}(\gamma_{j}) \gamma_{j} \nu_{j}.
\end{array}
\end{displaymath}

Take $\delta_{j}$ and $\delta_{j}^{'}$ to be the direct sum of $r_{j}$ and 
$s_{j}$ copies of $\gamma_{j}$, where $r_{j}$ and $s_{j}$ are to be specified.  
The condition, for $t = 0,1$, that
\begin{displaymath}
\delta_{j} \phi_{j}^{t} + \delta_{j}^{'} \phi_{j}^{1-t} = \phi_{j+1}^{t} \gamma_{j},
\end{displaymath}
understood up to unitary equivalence, then becomes the condition
\begin{displaymath}
r_{j} \gamma_{j} (l_{j}^{t} \mu_{j} + (k_{j} - l_{j}^{t}) \nu_{j}) + s_{j} 
\gamma_{j} (l_{j}^{1-t} \mu_{j} + (k_{j} - l_{j}^{1-t}) \nu_{j}) =
(l_{j+1}^{t} \mu_{j+1} + (k_{j+1} - l_{j+1}^{t}) \nu_{j+1}) \gamma_{j},
\end{displaymath}
also up to unitary equivalence.  As $    \mathrm{K}_{0}(\mu_{j})$ and $    \mathrm{K}_{0}(\nu_{j})$ are independent this is equivalent to the two equations
\begin{displaymath}
\begin{array}{rl}
r_{j} l_{j}^{t} + s_{j} l_{j}^{1-t} & = \mathrm{mult}(\gamma_{j}) l_{j+1}^{t}, \\
(r_{j} + s_{j}) k_{j} & = \mathrm{mult}(\gamma_{j}) k_{j+1}.
\end{array}
\end{displaymath}
Choose $r_{j} = 2 \mathrm{mult}(\gamma_{j})$ and $s_{j} = \mathrm{mult}(\gamma_{j})$, so that
\begin{displaymath}
k_{j+1} = 3 k_{j}
\end{displaymath}
and 
\begin{displaymath}
l_{j+1}^{t} = 2 l_{j}^{t} + l_{j}^{1-t}
\end{displaymath}

Taking $k_{1} = 1, l_{1}^{0} = 0,$ and $l_{1}^{1} = 1$ we have $k_{j} = 3^{j-1}$ 
and $ l_{j}^{1} - l_{j}^{0} = 1$ for all $j$ and, in particular, these quantities are non-zero, as 
required above.

Next let us show that, up to unitary equivalence preserving the equations 
\begin{displaymath}
\delta_{j} \phi_{j}^{t} + \delta_{j}^{'} \phi_{j}^{1-t} =
\phi_{j+1}^{t} \gamma_{j},
\end{displaymath}
$\phi_{j+1}^{0} \beta_{j}$
is a direct summand of $\delta_{j}^{'} = \mathrm{mult}(\gamma_{j}) \gamma_{j}$, 
and $\phi_{j+1}^{1} \beta_{j}$ is a direct summand of $\delta_{j} = 2 \mathrm{mult}(\gamma_{j}) \gamma_{j}$.

Note that $\phi_{j+1}^{t} \beta_{j}$ is the direct sum of $l_{j+1}^{t}$ copies 
of $p_{j+1} \otimes \beta_{j}$ and $(k_{j+1} - l_{j+1}^{t}) \mathrm{dim} (p_{j+1})$ copies of $\beta_{j}$, 
whereas $\delta_{j}^{'}$ and $\delta_{j}$ contain, respectively, $q_{j} \mathrm{mult} 
(\gamma_{j})$ and $2 q_{j} \mathrm{mult} (\gamma_{j})$ copies of $\beta_{j}$.  Note also that by Theorem 
8.1.5 of [Hu] that a trivial projection of dimension at least $\mathrm{dim} (p_{j+1}) 
+ \mathrm{dim} X_{j+1}$ in $\mathrm{C}(X_{j+1}) \otimes K$ contains a copy of $p_{j+1}$.  Therefore, 
$\mathrm{dim} (p_{j+1}) + \mathrm{dim} X_{j+1}$ copies of $\beta_{j}$ contain a copy 
of $p_{j+1} \otimes \beta_{j}$.  It follows that $k_{j+1} (2 \mathrm{dim} (p_{j+1}) + \mathrm{dim} X_{j+1})$
copies of $\beta_{j}$ contain a copy of $\phi_{j+1}^{t}$ when $t$ is either 1 or 0.  
Here, by a copy of a given map from $D_{j}$ to $D_{j+1}$ we mean another map obtained from it by conjugating
by a partial isometry in $D_{j+1}$ with initial projection the image of the unit.

Note that
\begin{displaymath}
\begin{array}{rl}
k_{j+1} (2 \mathrm{dim} (p_{j+1}) + \mathrm{dim} X_{j+1}) & = 3 k_{j} ( 2 \mathrm{mult}(\gamma_{j}) \mathrm{dim} (p_{j}) + n_{j} 
                                                        \mathrm{dim} X_{j}) \\
                                                 & \leq 3 k_{j} (2 \mathrm{dim} (p_{j}) + \mathrm{dim} X_{j}) \mathrm{mult}(\gamma_{j}),
\end{array}
\end{displaymath}
and that $k_{j}$, $\mathrm{dim} (p_{j})$ and $\mathrm{dim} X_{j}$ have already 
been specified and do not depend on $n_{j}$.  It follows that, with
\begin{displaymath}
q_{j} = 3 k_{j} ( 2 \mathrm{dim} (p_{j}) + \mathrm{dim} X_{j}),
\end{displaymath}
$q_{j} \mathrm{mult}(\gamma_{j})$ copies of $\beta_{j}$ contain a copy 
of $\phi_{j+1}^{t} \beta_{j}$ for $t = 0,1$.  In particular $\delta_{j}^{'}$ 
and $\delta_{j}$ contain copies, respectively,
of $\phi_{j+1}^{0} \beta_{j}$ and $\phi_{j+1}^{1} \beta_{j}$.

With this choice of $q_{j}$, let us show that for each $t = 0,1$ there exists a 
unitary $u_{t} \in D_{j+1}$ commuting with the image of $\phi_{j+1}^{t} \gamma_{j}$, such that 
$(\mathrm{Ad} u_{0}) \phi_{j+1}^{0} \beta_{j}$ is a direct summand of $\delta_{j}^{'}$ 
and $(\mathrm{Ad} u_{1}) \phi_{j+1}^{1} \beta_{j}$ is a direct summand of $\delta_{j}$.
In other words, for each $t = 0,1$ we must show that the partial isometry constructed 
in the preceding paragraph, producing a copy of $\phi_{j+1}^{t} \beta_{j}$ inside either
$\delta_{j}^{'}$ or $\delta_{j}$, may be chosen in such a way that it extends to a 
unitary element of $D_{j+1}$ --- which in addition commutes with the image of $\phi_{j+1}^{t} \gamma_{j}$.

We will consider the case $t=0$.  The case $t=1$ is similar.  Let us first show that 
the partial isometry in $D_{j+1}$, transforming $\phi_{j+1}^{0} \beta_{j}$ into a direct
summand of $\delta_{j}^{'}$, may be chosen to lie in the commutant of the image of 
$\phi_{j+1}^{0} \gamma_{j}$.  Note first that the unit of the image of $\phi_{j+1}^{0} \beta_{j}$
--- the initial projection of the partial isometry --- lies in the commutant of the 
image of $\phi_{j+1}^{0} \gamma_{j}$.  Indeed, this projection is the image
by $\phi_{j}^{1}$ of the unit of $C_{j}$.  The property that $\beta_{j} \phi_{j}^{1}$ is 
a direct summand of $\gamma_{j}$ implies in particular that the image by
$\beta_{j} \phi_{j}^{1}$ of the unit of $C_{j}$ commutes with the image of $\gamma_{j}$.  
The image by $\phi_{j+1}^{0} \beta_{j} \phi_{j}^{1}$ of the unit of $C_{j}$
(i.e. the unit of the image of $\phi_{j+1}^{0} \beta_{j}$) therefore commutes with the 
image of $\phi_{j+1}^{0} \gamma_{j}$, as asserted.

Note also that the final projection of the partial isometry commutes with the 
image of $\phi_{j+1}^{0} \gamma_{j}$.  Indeed, it is the unit of the image of a direct summand
of $\delta_{j}^{'}$, and since $D_{j}$ is unital it is the image of the unit of 
$D_{j}$ by this direct summand;  since $C_{j}$ is unital and $\phi_{j}^{1}: C_{j} \longrightarrow D_{j}$
is unital, the projection in question is the image of the unit of $C_{j}$ by a direct 
summand of $\delta_{j}^{'} \phi_{j}^{1}$.  But $\delta_{j}^{'} \phi_{j}^{1}$ is itself a direct
summand of $\phi_{j+1}^{0} \gamma_{j}$ (as $\phi_{j+1}^{0} \gamma_{j} = \delta_{j} 
\phi_{j}^{0} + \delta_{j}^{'} \phi_{j}^{1}$), and so the projection in question is the image of the unit
of $C_{j}$ by a direct summand of $\phi_{j+1}^{0} \gamma_{j}$, and in particular 
commutes with the image of $\phi_{j+1}^{0} \gamma_{j}$.            

Note that both direct summands of $\phi_{j+1}^{0} \gamma_{j}$ under consideration 
($\phi_{j+1}^{0} \beta_{j} \phi_{j}^{1}$ and a copy of it) factor through the evaluation of $C_{j}$ 
at the point $x_{j}$, and so are contained in the largest such direct summand of 
$\phi_{j+1}^{0} \gamma_{j}$;  this largest direct summand, say $\pi_{j}$, is seen to
exist by inspection of the construction of $\phi_{j+1}^{0} \gamma_{j}$.  Since both 
projections under consideration (the images of the unit of $C_{j}$ by the two copies of 
$\phi_{j+1}^{0} \beta_{j} \phi_{j}^{1}$) are less than $\pi_{j}(1)$, to show that 
they are unitarily equivalent in the commutant of the image of $\phi_{j+1}^{0} \gamma_{j}$ 
(in $D_{j+1}$) it is sufficient to show that they are unitarily equivalent in the 
commutant of the image of $\pi_{j}$ in $\pi_{j}(1) D_{j+1} \pi_{j}(1)$.  Note that this image
is isomorphic to $\mathrm{M}_{\mathrm{dim} p_{j}}(C)$.  By construction, the two projections 
in question are Murray-von Neumann equivalent --- in $D_{j+1}$ and therefore
in $\pi_{j}(1) D_{j+1} \pi_{j}(1)$ --- but all we shall use from this is that 
they have the same class in $K^{0} X_{j+1}$.  Note that the dimension of these projections is 
$(k_{j+1} \mathrm{dim} (p_{j+1}))(k_{j} \mathrm{dim} (p_{j}))$, and that the dimension 
of $\pi_{j}(1)$ is $k_{j+1} \mathrm{dim} (p_{j+1}) + l_{j+1}^{0} {(\mathrm{dim} (p_{j+1}))}^{2}$.
Since the two projections under consideration commute with $\pi_{j}(C_{j})$, and this 
is isomorphic to $\mathrm{M}_{\mathrm{dim} (p_{j})}(C)$, to prove unitary equivalence in the commutant of 
$\pi_{j}(C_{j})$ in $\pi_{j}(1) D_{j+1} \pi_{j}(1)$ it is sufficient to prove unitary 
equivalence of the product of these projections with a fixed minimal projection of $\pi_{j}(C_{j})$,
say $e$.  Since $K^{0} X_{j+1}$ is torsion free, the products of the two projections 
under consideration with $e$ still have the same class in $K^{0} X_{j+1}$.  To prove 
that they are unitarily equivalent in $e D_{j+1} e$ , it is sufficient (and necessary) to prove 
that both they and their complements inside $e$ are Murray von-Neumann equivalent.  Since both 
the cut-down projections and their complements inside $e$ have the same class in $K^{0} X_{j+1}$, to 
prove that the two pairs are equivalent it is sufficient, by Theorem 8.1.5 of [Hu], to show that all four projections 
have dimension at least $\frac{1}{2} \mathrm{dim} X_{j+1}$.  Dividing the dimensions 
above by $\mathrm{dim} (p_{j})$ (the order of the matrix algebra), we see that the dimension
of the first pair of projections is $ k_{j+1} k_{j} \mathrm{dim} (p_{j+1}) = k_{j+1} k_{j} 
\mathrm{mult}(\gamma_{j}) \mathrm{dim} (p_{j})$.  The dimension of $e$ is 
$k_{j+1} \mathrm{mult}(\gamma_{j}) +l_{j+1}^{0} \mathrm{mult}(\gamma_{j}) \mathrm{dim} (p_{j+1})$, 
so that the dimension of the second pair of projections is
$\mathrm{mult}(\gamma_{j}) (k_{j+1} + l_{j+1}^{0} \mathrm{dim} (p_{j+1}) - k_{j+1} k_{j} \mathrm{dim} 
(p_{j}))$.  Since $\mathrm{dim} (p_{1}) \geq \frac{1}{2} \mathrm{dim} X_{1}$, 
$\mathrm{dim} (p_{j+1}) = \mathrm{mult}(\gamma_{j}) \mathrm{dim} (p_{j})$, $\mathrm{dim} 
X_{j+1} = n_{j} \mathrm{dim} X_{j}$, and $\mathrm{mult}(\gamma_{j}) \geq n_{j}$ (for all $j$),
we have $\mathrm{dim} (p_{j+1}) \geq \frac{1}{2} \mathrm{dim} X_{j+1}$ (for all $j$).  
Since $k_{j+1} k_{j}$ is non-zero for all $j$, the first inequality holds.  Since $l_{j+1}^{0}$ 
is non-zero for all $j$, the second inequality holds if $\mathrm{mult}(\gamma_{j})$ is 
strictly greater than $k_{j+1} k_{j}$.  (One then has, using $\mathrm{dim} (p_{j+1}) =
\mathrm{mult}(\gamma_{j}) \mathrm{dim} (p_{j})$ twice, that the dimension of the 
second pair of projections is at least $\mathrm{dim} (p_{j+1})$.)  Since $k_{j+1} k_{j} = 3 {k_{j}}^{2}$, 
and $k_{j}$ was specified before $n_{j}$, we may modify the choice of $n_{j}$ so 
that $\mathrm{mult}(\gamma_{j})$ --- which is greater than $n_{j}$ --- is sufficiently large.

This shows that the two projections in $D_{j+1}$ under consideration are unitarily equivalent 
by a unitary in the commutant of the image of $\phi_{j+1}^{0} \gamma_{j}$.  
Replacing $\phi_{j+1}^{0}$ by its composition with the corresponding inner automorphism, we 
may suppose that the two projections in question are equal.  In other words
$\phi_{j+1}^{0} \beta_{j}$ is unitarily equivalent to the cut-down of $\delta_{j}^{'}$ by 
the projection $\phi_{j+1}^{0} \beta_{j}(1)$.

Now consider the compositions of these two maps with $\phi_{j}^{1}$, namely $\phi_{j+1}^{0} 
\beta_{j} \phi_{j}^{1}$ and the cut-down of $\delta_{j}^{'} \phi_{j}^{1}$
by the projection $\phi_{j+1}^{0} \beta_{j}(1)$.  Since both of these maps can be viewed as 
the cut-down of $\phi_{j+1} \gamma_{j}$ by the same projection, they are in fact
the same map.  Thus any unitary inside the cut-down of $D_{j+1}$ by $\phi_{j+1}^{0} \beta_{j}(1)$ 
taking $\phi_{j+1}^{0} \beta_{j}$ into the cut-down of $\delta_{j}^{'}$ by
this projection (such a unitary is known to exist) must commute with the image of 
$\phi_{j+1}^{0} \beta_{j} \phi_{j}^{1}$ and hence with the image of $\phi_{j+1}^{0} \gamma_{j}$,
since this commutes with the projection $\phi_{j+1}^{0} \beta_{j}(1) = \phi_{j+1}^{0}(\beta_{j} 
\phi_{j}^{1}(1))$.  The extension of such a partial unitary to a unitary $u_{0}$
in $D_{j+1}$ equal to one inside the complement of this projection then belongs to the 
commutant of the image of $\phi_{j+1}^{0} \gamma_{j}$, and transforms $\phi_{j+1}^{0} \beta_{j}$
into the cut-down of $\delta_{j}^{'}$ by this projection, as desired.

As stated above, the proof for the case $t=1$ is similar.

Inspection of the construction of the maps $\delta_{j}^{'} - \phi_{j}^{0} \beta_{j}$ and 
$ \delta_{j} - \phi_{j}^{1} \beta_{j}$ shows that they are injective, as required by the hypotheses of 
section 2.

Replacing $\phi_{j+1}^{t}$ with $(\mathrm{Ad} u_{t}) \phi_{j+1}^{t}$,
we have an inductive sequence
\begin{displaymath}
A_{1} \stackrel{\theta_{1}}{\longrightarrow} A_{2} \stackrel{\theta_{2}}{\longrightarrow} \cdots
\end{displaymath}
satisfying the hypotheses of section 2.  (The existence of
$\alpha_{j}$ homotopic to $\beta_{j}$ and 
non-zero on a given element of $D_{j}$, defined by another 
point evaluation, is clear.)

By Theorem 2.3 there exists a sequence 
\begin{displaymath}
A_{1} \stackrel{{\theta}_{1}^{'}}{\longrightarrow} A_{2} \stackrel{{\theta}_{2}^{'}}{\longrightarrow} \cdots ,
\end{displaymath}
with $\theta_{j}^{'}$ homotopic to $\theta_{j}$ (and so agreeing with $\theta_{j}$ on $\mathrm{K}_{0}$), 
the inductive limit of which is simple.

Since the map $    \mathrm{K}_{0}(\theta_{j}^{'})$ (considered as a map between single copies of the 
integers) takes the canonical generator $1 \in \Z$ to $L_{j+1}$, we may conclude that the
simple inductive limit in question has the desired $    \mathrm{K}_{0}$-group.  That the positive 
elements are all those greater than $k$ follows from the fact that at each stage, $l_{j} + 1$ is
the smallest positive element in $    \mathrm{K}_{0} A_{j} = \Z$ and
\begin{displaymath}
\mathrm{lim} \frac{l_{j} + 1}{\Pi_{k=1}^{j} L_{k}} = 
\mathrm{lim} \frac{a \Pi_{k=2}^{j} L_{j} +1}{b \Pi_{k=2}^{j} L_{j}}
= k + \mathrm{lim} \frac{1}{\Pi_{k=1}^{j} L_{k}} = k.
\end{displaymath}
Theorem 3.1 follows. \hfill $\Box$

Finally, one might reasonably ask whether ${\mathrm{K}_{0}(A_{(\mathbf{n},k)})}^{+}$ can be made to
contain $k$.  There is no reason \emph{a priori} why this should not be possible, but the construction 
above does not seem amenable to modifications which would achieve this result.  Roughly speaking, the 
$\mathrm{K}_0$-group in Theorem 3.1 can be thought of as an inductive limit of sub-ordered groups of
ordered $\mathrm{K}_0$-groups of homogeneous $C^*$-algebras.  In order that the inductive limit of
Theorem 3.1 be simple, one must introduce point evaluations via the maps $\beta_j$.  In the absence 
of these point evaluations, one could have maps $\Psi: \mathsf{Z}_{mk} \rightarrow \mathsf{Z}_{mnk}$ 
with $\Psi(nk)=mnk$ at the level of $\mathrm{K}_0$ between the building blocks $A_i$ and $A_{i+1}$.  
With these point evaluations, however, one is forced into a situation where $\Psi(nk)$ is necessarily
strictly less than $mnk$.

\end{document}